\documentclass[11pt]{elsarticle}
\usepackage[utf8]{inputenc}
\usepackage[margin=3cm]{geometry}
\usepackage{amsmath,amssymb,amsthm,bbm}
\usepackage[table]{xcolor}
\usepackage{color,soul}
\usepackage{graphicx}
\usepackage{booktabs}
\usepackage{todonotes}
\usepackage{float}
\usepackage{latexsym}
\usepackage{bm}
\usepackage{multirow}
\usepackage[shortlabels]{enumitem}
\usepackage{overpic}

\journal{Applied Mathematics Letters}

\begin{document}

\begin{frontmatter}

\title{Single-spike solutions to the 1D shadow Gierer-Meinhardt problem}

\author[label1]{Annalisa Iuorio}\ead{annalisa.iuorio@univie.ac.at}
\author[label2]{Christian Kuehn}\ead{ckuehn@ma.tum.de}

\affiliation[label1]{organization={University of Vienna, Faculty of Mathematics},
            addressline={Oskar-Morgenstern-Platz 1}, 
            city={Vienna},
            postcode={1090},
            country={Austria}}
            
\affiliation[label2]{organization={Technical University of Munich},
            addressline={Boltzmannstr. 3}, 
            city={Garching bei M\"unchen},
            postcode={85748},
            country={Germany}}

\begin{abstract}
%WHY: open relevant problem still no explicit formulas
%HOW: using generalized hyperbolic functions
%WHAT: we give explicit formulas for single spike solutions
%SO WHAT: fills a gap and provides inspiration/guide for spike solutions in multiD

A fundamental example of reaction-diffusion system exhibiting Turing type pattern formation is the Gierer-Meinhardt system, which reduces to the shadow Gierer-Meinhardt problem in a suitable singular limit. Thanks to its applicability in a large range of biological applications, this singularly perturbed problem has been widely studied in the last few decades via rigorous, asymptotic, and numerical methods. However, standard matched asymptotics methods do not apply \cite{Ni_1998,Wei_1998_2}, and therefore  analytical expressions for single spike solutions are generally lacking.\\ 
By introducing an ansatz based on generalized hyperbolic functions, we determine exact radially symmetric solutions to the one-dimensional shadow Gierer-Meinhardt problem for any $1 < p < \infty$, representing both inner and boundary spike solutions depending on the location of the peak. Our approach not only confirms numerical results existing in literature, but also provides guidance for tackling extensions of the shadow Gierer-Meinhardt problem based on different boundary conditions (e.g.~mixed) and/or $n$-dimensional domains.
\end{abstract}

%%%Graphical abstract
%\begin{graphicalabstract}
%\textcolor{red}{TODO?}
%\end{graphicalabstract}

%%%Research highlights
%\begin{highlights}
%\item \textcolor{red}{TODO}
%\item \textcolor{red}{TODO}
%\end{highlights}

\begin{keyword}
Singular perturbation, explicit solution, Gierer-Meinhardt, nonlinear elliptic PDE, boundary value problem. 
\vspace{.2cm}
%\MSC[2020] 34B16 \sep 34E13 \sep 34E15 \sep 35B25 \sep 35B40 \sep 65L11 \sep 76A30.

\end{keyword}

\end{frontmatter}

\section{Introduction}

In this paper, we provide an analytical expression for single spike solutions to the shadow Gierer-Meinhardt equation on a one-dimensional domain $\Omega = [-L,L] \subset \mathbb{R}$ (with $L>0$), given by
\begin{equation} \label{eq:Ni}
\begin{aligned}
    &\varepsilon^2 \Delta u - u +u^p = 0 \quad \,\,\, \text{ in } \Omega, \\
    &u > 0 \quad \text{ in } \Omega \quad \text{ and } \quad \frac{\partial u}{\partial \nu} = 0 \quad \text{ on } \partial \Omega,
    %, \\
    %\frac{\partial u}{\partial \nu} &= 0, \quad %\text{ on }\partial\Omega,
\end{aligned}
\end{equation}
where $0 < \varepsilon \ll 1$ and $1<p< \infty$ (see \cite{Wei_2008} and references therein). 

This problem arises from the full Gierer-Meinhardt (GM) model, a fundamental example of a Turing type reaction-diffusion system describing the dynamics of an activator $a(x,t)$ and an inhibitor $h(x,t)$ \cite{Gierer_1972}. As these two variables (are assumed to) evolve on different time scales (\emph{slowly} diffusing activator, \emph{fast} diffusing inhibitor), the parameter $D$ (representing the ratio between $h$ and $a$ diffusion coefficients, respectively) satisfies $D \gg \varepsilon$ (see \cite{Wei_2008}). Therefore, it is possible to study the GM model in the singular limit $D \to \infty$ by rescaling both $a$ and $h$ via a function $u$, which leads to the shadow Gierer-Meinhardt system in \eqref{eq:Ni}.

%(see e.g.~\cite{Cerami_1998,Gui_1999,Kowalczyk_1999,Wei_1998_3,Wei_1998_2} for inner spike solutions, \cite{del_Pino_1999,Gui_1996,Gui_2000,Ni_1995_2,Wei_1997,Wei_1999_2,Yanyan_1998} for boundary spike solutions

The existence and main properties of such solutions for general $n$-dimensional domains has been extensively studied in the last few decades (see e.g.~\cite{Gui_1999,Kowalczyk_1999,Wei_1998_2} for inner spike solutions, \cite{Gui_2000,Ni_1995_2,Wei_1999_2,Yanyan_1998} for boundary spike solutions, and \cite{Iron_2000,WEI_1999} for stability of spike solutions). Moreover, numerical profiles have been computed for different $p$-values using suitable boundary value problem solvers (see \cite{Ward_1996}). However, an explicit expression for inner spike solutions to \eqref{eq:Ni} is only available for $p=2$, where (see e.g. \cite{Iron_2000})
\begin{equation} \label{eq:p2}
u(x)=\frac32\, \mathrm{sech}^2\left( \frac{x}{2\, \varepsilon} \right).
\end{equation}
Using the method based on generalized hyperbolic functions introduced in \cite{Pandir_2013}, we can determine an analytical expression for inner spike solutions to \eqref{eq:Ni} for any $1<p<\infty$.

%\begin{description}
% \item[Boundary layers.] Generally occur with Neumann boundary conditions, but can occur also with Dirichlet BC. Moreover, the order of boundary spikes is \emph{algebraic} (see \cite{Gui_2000_2}).
% \item[Interior layers.] Can only occur with Dirichlet boundary conditions (see \cite{Kowalczyk_1999}). Moreover, the order of interior spikes is \emph{exponentially small} (see \cite{Gui_2000_2}).
%\end{description}

\section{Spike solutions}

Introducing the rescaled variable $\rho=\frac{x}{\varepsilon}$, a spike solution to \eqref{eq:Ni} is a radially symmetric monotone function $u(\rho)$ satisfying
\begin{subequations} \label{eq:inspN1}
\begin{align}
    &u''-u+u^p=0, \label{eq:inspN1eq}\\
    &u'(\rho^\ast)=0, \label{eq:inspN1peak}\\
    &u(\rho) \rightarrow 0 \text{ as } \rho \rightarrow \rho_l, \label{eq:inspN1limit}\\
    &u > 0 \label{eq:inspN1pos},
\end{align}    
\end{subequations}
where $\rho^\ast$ is the $\rho$ value where the peak is located and $\rho_l$ represents the limiting $\rho$-value where $u$ vanishes. In particular, for inner spike solutions we have $\rho^\ast=0$, $\rho_l=\infty$, whereas for boundary spike solutions $\rho^\ast=L$ and $\rho_l=0$ as $L \to \infty$. We note that on a general $n$-dimensional domain an extra term $\frac{n-1}{\rho}\,u'$ would appear on the left hand-side of Equation \eqref{eq:inspN1eq}, which drops out here as we consider $n=1$ (see e.g.~\cite{Ward_1996}). Therefore, an extension of the technique presented here would need to be considered in order to study Equation \eqref{eq:inspN1} on higher dimensional domains ($n>1$).

In order to determine an analytic expression for the solution to \eqref{eq:inspN1}, we introduce the ansatz
\begin{equation} \label{eq:1pansatz}
    u(\rho)=\frac{A}{\mathrm{cosh}_a^s(\rho)},
\end{equation}
where $\mathrm{cosh}_a(\rho)=\frac{m a^{k\rho}+q a^{-k\rho}}{2}$, and $A$, $a$, $m$, $q$, $k$, $s$ are constants to be derived (see \cite{Pandir_2013}). As long as the quantities $m$, $k\, \mathrm{log}\,a$, and $s$ are positive, such function decays at infinity, and is therefore a valid candidate for our problem.\\
\noindent
Plugging \eqref{eq:1pansatz} into \eqref{eq:inspN1}, we get
\begin{equation}
    A\,k^2\,s\,\mathrm{ln}^2(a)\,\mathrm{cosh}_a^{-2-s}(\rho) \left( s\,\mathrm{cosh}_a^2(\rho)-m\,q\left( 1+s \right) \right)-\frac{A}{\mathrm{cosh}_a^s(\rho)}+\frac{A^p}{\mathrm{cosh}_a^{p s}(\rho)}=0.
\end{equation}
By suitably rearranging terms, we obtain
\begin{equation} \label{eq:1phalf}
    A\,\mathrm{cosh}_a^{-s}(\rho)\left( -1+A^{p-1}\mathrm{cosh}_a^{-s(p-1)}(\rho)+k^2\,s\,\mathrm{ln}^2(a)\, \mathrm{cosh}_a^{-2}(\rho) \left( s\,\mathrm{cosh}_a^2(\rho)-m\,q\left( 1+s \right) \right) \right) = 0.
\end{equation}
In order to balance terms, a convenient choice is $s=\frac{2}{p-1}$, $a=e$. This reduces Equation \eqref{eq:1phalf} to
\begin{equation}
    A\,\mathrm{cosh}_a^{2}(\rho)\,\left( 4k^2-(p-1)^2 \right)+\left( A^p(p-1)^2-2A\,k^2\,m\,q(1+p) \right)=0.
\end{equation}
Since the first term depends on $\rho$ while the second one does not, we get
\begin{equation}
    \begin{aligned}
        4k^2-(p-1)^2 &=0,\\
        A^p(p-1)^2-2A\,k^2\,m\,q(1+p) &= 0,
    \end{aligned}
\end{equation}
which in turn leads to
\begin{equation} \label{eq:km}
    k=\frac12(p-1), \qquad m=\frac{2 A^{p-1}}{q\,(p+1)}.
\end{equation}
We notice that, according to the fixed parameters, the expressions for $m$, $k\, \mathrm{log}\,a$, and $s$ satisfy the sign constraints which ensure the decay at infinity of our candidate solution \eqref{eq:1pansatz}.\\
The last constraint we need to impose from \eqref{eq:inspN1} is $u'(\rho^\ast)=0$. Taking into account Equation \eqref{eq:inspN1pos}, this leads to
\begin{equation} \label{eq:A}
    A=e^{-\rho^\ast} \left( \frac{2}{q^2(1+p)} \right)^{\frac{1}{1-p}}.
\end{equation}
When we plug \eqref{eq:A} into \eqref{eq:km}, we are left with $m=e^{\rho^\ast(1-p)} q$. 
Substituting all these findings into \eqref{eq:1pansatz}, we obtain that
\begin{equation} \label{eq:1psol}
    u(\rho)=\left( \frac{1+\mathrm{cosh}\left( (1-p) (\rho-\rho^\ast) \right)}{1+p} \right)^{\frac{1}{1-p}}=:u_s(\rho),
\end{equation}
where the $q$ term cancels out. Hence, multiple choices of $q$ (and therefore different $m$ and $A$, correspondingly) lead to the same result, as long as the constraints in \eqref{eq:km} and \eqref{eq:A} hold. Equation \eqref{eq:inspN1limit} is satisfied both for inner ($\rho_l=\infty$) and boundary ($\rho_l=0$ as $L \to \infty$) spike solutions.\\
Rescaling the independent variable to the original one, we have that
\begin{equation} \label{eq:1psolx}
    u(x)=\left( \frac{1+\mathrm{cosh}\left( (1-p) \left( \frac{|x-x^\ast|}{\varepsilon} \right) \right)}{1+p} \right)^{\frac{1}{1-p}}=:u_s(x),
\end{equation}
where $x^\ast=\rho^\ast \varepsilon$. We observe that for $p=2$ this expression coincides with the one in \eqref{eq:p2}.

\section{Numerical solutions} \label{ssec:num}

\begin{figure}[!h]
    \centering
    \begin{overpic}[scale=0.55]{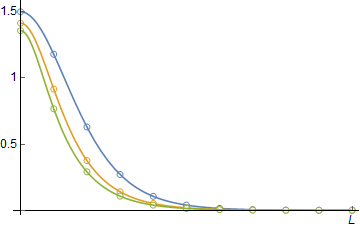}
	\put(-5,32){$u$}
	\put(50,-2){$\rho$}
	\put(-9,65){(a)}
	\end{overpic}
	\hspace{.2cm}
	\begin{overpic}[scale=0.55]{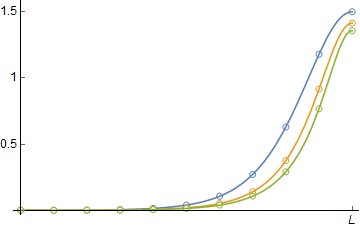}
	\put(-5,32){$u$}
	\put(50,-2){$\rho$}
	\put(-9,65){(b)}
	\end{overpic}
    \caption{Comparison between the analytical solution to \eqref{eq:inspN1} obtained in \eqref{eq:1psol} (continuous line) and the numerical one (circles) obtained with the shooting strategy described in Section \ref{ssec:num} for in the case $p=2$ (blue), $p=3$ (orange), and $p=4$ (green). The panels (a), (b) display inner and  boundary spike solutions, respectively.}
    \label{fig:plot1pcomp}
\end{figure}
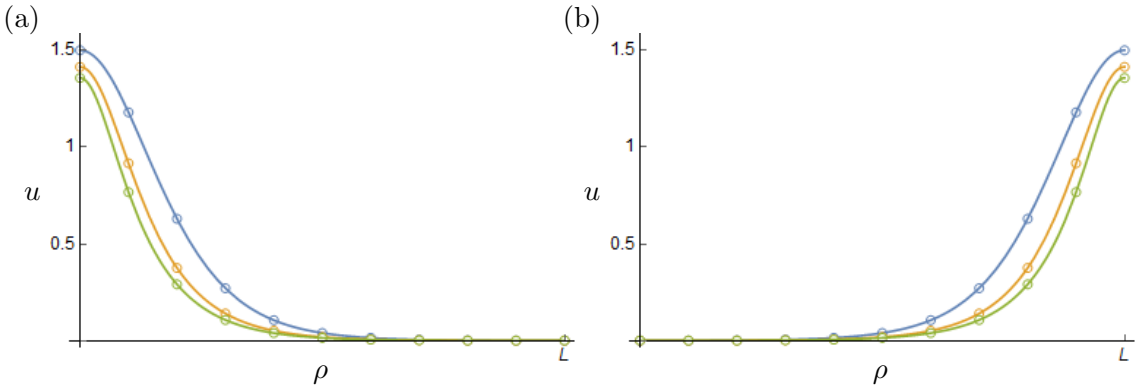

The accuracy of our solution for inner spike solutions in the case of a generic $1<p<\infty$ is confirmed by a qualitative comparison with numerical solutions obtained using the NDSolve routine in Mathematica. For such problems, this function is based on a combination of the LSODA algorithm and Newton's method.\\
Due to the highly singular nature of the problem, our numerical approach is based on a shooting strategy as follows:
\begin{description}
\item[Step 1.] We reformulate Equation \eqref{eq:inspN1eq} as a system of two first-order ODEs:
\begin{equation} \label{eq:inspSYS}
\begin{aligned}
    u' &= v,\\
    v' &= u-u^p.
\end{aligned}
\end{equation}
Moreover, we equip Equation \eqref{eq:inspSYS} with the boundary conditions
\begin{subequations} \label{eq:SYSBC}
\begin{align}
    u(\rho^\ast) &= a, \label{eq:SYSBCu}\\
    v(\rho^\ast) &= 0, \label{eq:SYSBCv}
\end{align}    
\end{subequations}
where Equation \eqref{eq:SYSBCv} corresponds to Equation \eqref{eq:inspN1peak} and $a$ lies in an interval of radius $\delta$ centered at the analytical value $u_s(\rho^\ast)$, i.e. $a \in [u_s(\rho^\ast)-\delta,\,u_s(\rho^\ast)+\delta]$, with $\delta = 0.1$.
\item[Step 2.] We solve Equation \eqref{eq:inspSYS}-\eqref{eq:SYSBC} using the NDSolve routine in Mathematica and pick the (unique) solution which verifies the artificial boundary condition (see \cite{Ward_1996})
\begin{equation} \label{eq:artBC}
    u(\rho_l)+v(\rho_l) \leq \eta,
\end{equation}
where we fix $\eta=0.01$.
\item[Step 3.] We extract the corresponding value of $a$ which satisfies Equation \eqref{eq:artBC}.
\end{description}
%the numerical results obtained in \cite{Iron_2000}
A comparison between the analytical expression in \eqref{eq:1psol} and the numerical solution obtained with the above strategy for $p=2,\,3,\,4$ is illustrated in Figure \ref{fig:plot1pcomp}.
\vspace{.5cm}

\section{Conclusion}

In this work, we provide an explicit formulation for single-spike solutions to the 1D shadow Gierer-Meinhardt problem for any $1<p<\infty$. This formula (which returns the known expression in the case $p=2$ and matches the numerical findings both for inner and boundary spikes) constitutes an invaluable tool to improve our understanding of such a crucial problem describing pattern formation. This improved knowledge allows us not only to tackle extensions of the shadow Gierer-Meinhardt problem (e.g.~different boundary conditions, higher-dimensional domains), but also to shed a light on other relevant problems with similar features. 

\section{Acknowledgements}
AI acknowledges support from an FWF Hertha Firnberg Research Fellowship (T 1199-N).

\bibliographystyle{plain}
\bibliography{references}

\end{document}